\documentclass[12pt]{amsart}
\usepackage{amscd,indentfirst,amsmath,amscd,latexsym,amsthm,amsfonts,amssymb,graphicx,color,geometry,hyperref}
\usepackage{mathtools}
\usepackage[utf8]{inputenc}
\usepackage{amsmath}
\usepackage{hyperref}
\usepackage{enumerate}
\usepackage{comment}
\usepackage{amsfonts}
\usepackage{hyperref}
\usepackage{amssymb}
\usepackage{color}
\usepackage{graphicx}
\usepackage{verbatim}

\newcommand{\1}{\textbf{1}}
\newcommand{\Lmn}{{\mathcal{L}(^m\R^n)}}

\newcommand{\ext}{\text{ext}}

\newcommand{\R}{\mathbb{R}}

\newcommand{\mf}{\mathfrak{m}}
\newcommand{\bb}{\begin{equation}}
\newcommand{\ee}{\end{equation}}

\begin{document}
\title[The optimal Bohnenblust-Hille constants: a computational solution]{The optimal multilinear Bohnenblust-Hille constants:\\
a computational solution for the real case}
\author[Fernando Vieira Costa J\'unior]{Fernando Vieira Costa J\'unior}
\address{Departamento de Matem\'{a}tica,
	Universidade Federal da Para\'{\i}ba,
	58.051-900 - Jo\~{a}o Pessoa, Brazil.}
\email{costa.junior.ufpb@gmail.com}
\address{Universit\'e Clermont Auvergne,
	CNRS, LMBP, F-63000 Clermont-Ferrand, France.}
\email{fernando.vieira\underline{ }costa\underline{ }junior@etu.uca.fr}
\thanks{2010 Mathematics Subject Classification: 15-04; 46-04}
\thanks{The author is supported by Capes}
\keywords{Bohnenblust-Hille inequality, optimal constants}

\begin{abstract}
The Bohnenblust-Hille inequality for $m$-linear forms was proven in 1931 as
a generalization of the famous 4/3-Littlewood inequality. The optimal
constants (or at least their asymptotic behavior as $m$ grows) is unknown,
but significant for applications. A recent result, obtained by Cavalcante, Teixeira and
Pellegrino, provides a kind of algorithm, composed by finitely many
elementary steps, giving as the final outcome the optimal truncated
Bohnenblust-Hille constants of any order. But the procedure of Cavalcante
\textit{et al.} has a fairly large number of calculations and computer
assistance cannot be avoided. In this paper we present a computational
solution to the problem, using the Wolfram Language. We also use this
approach to investigate a conjecture raised by Pellegrino and Teixeira,
asserting that $C_m=2^{1-1/m}$ for all $m\in\mathbb{N}$ and to reveal
interesting unknown facts about the geometry of $B_{\mathcal{L}(^3\mathbb{R}^3)}$.
\end{abstract}

\maketitle

\section{Introduction and notation}

The Bohnenblust-Hille inequality was proven in 1931 (see \cite{bohnenblust1931absolute}) by H. F. Bohnenblust (1906-2000) and C. E. Hille
(1894-1980) as a generalization of the famous 4/3-Littlewood inequality. For
the real case, it states that for each $m\in \mathbb{N}$ there exists an
optimal constant $C_{m}\geq 1$ such that, for all $n\in \mathbb{N}$ and all $m$-linear forms $T:{\mathbb{R}^{n}\times \cdots \times \mathbb{R}^{n}\rightarrow \mathbb{R}}$, we have
\begin{equation}
\biggl(\sum_{i_{1},...,i_{m}}^{n}|T(e_{i_{1}},...,e_{i_{m}})|^{\frac{2m}{m+1}}\biggr)^{\frac{m+1}{2m}}\leq C_{m}\sup \{|T(x_{1},...,x_{n})|\in \mathbb{R}:\Vert x_{i}\Vert \leq 1\}.  \label{bhinequality}
\end{equation}

The numbers $C_{m}$, $m\in \mathbb{N}$, are called \emph{Bohnenblust-Hille
constants}. If we fix $n$ in the assumptions above, then we have the \emph{truncated Bohnenblust-Hille constants}, denoted by $C_{m}^{(n)}$, which
satisfy
\begin{equation}
\biggl(\sum_{i_{1},...,i_{m}=1}^{n}|T(e_{i_{1}},...,e_{i_{m}})|^{\frac{2m}{m+1}}\biggr)^{\frac{m+1}{2m}}\leq C_{m}^{(n)}\sup \{|T(x_{1},...,x_{n})|\in\mathbb{R}:\Vert x_{i}\Vert \leq 1\}  \label{bhtruncatedinequality}
\end{equation}
for all $m$-linear forms $T:{\mathbb{R}^{n}\times \cdots \times \mathbb{R}^{n}\rightarrow \mathbb{R}}$. In this case, for each $m\in \mathbb{N}$, we
have $C_{m}=\inf_{n}C_{m}^{(n)}$. For applications of the Bohnenblust-Hille inequality in physics we refer to
\cite{montanaro2012some} and references therein.

Let us establish some notations that will be carried out along this work.
The transpose of a matrix $M$ will be denoted by $M^{t}$. For $k\in \mathbb{N}$, the $k$-th cartesian product of a set $A$ will be denoted by $A^{k}$.
For example, $[\mathbb{R}^{n}]^{m}=\mathbb{R}^{n}\times \overset{m}{\cdots }\times \mathbb{R}^{n}$, $m,n\in \mathbb{N}$. If $E$ is a normed space, $B_{E}$ and $\text{ext}(B_{E})$ will denote the closed unit ball and the set of
its extreme points, respectively. Given $m,n\in \mathbb{N}$, the $n$-dimensional Euclidean space $\mathbb{R}^{n}$ will be endowed with the norm
\begin{equation*}
\Vert x\Vert =\sup \{\Vert x^{(i)}\Vert \in \mathbb{R}:i=1,...,n\},~\text{for all}~x=(x^{(1)},...,x^{(n)})\in \mathbb{R}^{n},
\end{equation*}
and the space ${\mathcal{L}(^{m}\mathbb{R}^{n})}$ of all $m$-linear forms $T:[\mathbb{R}^{n}]^{m}\rightarrow \mathbb{R}$ will be equipped with the norm
\begin{equation*}
\Vert T\Vert =\sup \{|T(x_{1},...,x_{n})|\in \mathbb{R}:x_{i}\in \mathbb{R}^{n}~\text{and}~\Vert x_{i}\Vert \leq 1,i=1,...,n\},~\text{for all}~T\in {\mathcal{L}(^{m}\mathbb{R}^{n})}.
\end{equation*}
From the Krein-Milman theorem we obtain the formula
\begin{equation}
\Vert T\Vert =\sup \{|T(x_{1},...,x_{n})|\in \mathbb{R}:x_{i}\in \text{ext}(B_{\mathbb{R}^{n}}),i=1,...,n\},~\text{for all}~T\in {\mathcal{L}(^{m}\mathbb{R}^{n})}.  \label{normTformula0}
\end{equation}
For a complete discussion about extreme points of convex sets and the
Krein-Milman (and Minkowski) theorem, we refer to \cite[Chapter 8]{simon2011convexity}. In the Example $8$ of the referenced chapter, we see that, for
each $n\in \mathbb{N}$, the extreme points of the closed unit ball $B_{\mathbb{R}^{n}}$ in the $n$-dimensional Euclidean space $\mathbb{R}^{n}$
equipped with the sup norm is simply
\begin{equation}
\text{ext}(B_{\mathbb{R}^{n}})=\{(\alpha _{1},...,\alpha _{n})\in \mathbb{R}^{n}:\alpha _{i}=1~\text{or}~\alpha _{i}=-1,~\text{for each}~i=1,...,n\}.
\label{extBRnm}
\end{equation}

\section{Bounds for the Bohnenblust-Hille constants}

The Bohnenblust--Hille inequalities are part of a bigger family of inequalites called Hardy--Littlewood inequalities. Through the last years, researchers have chase for better upper bounds for the constants in these inequalities (for example \cite{bayart2014bohr}, \cite{araujo2014upper} and \cite{pellegrino2017towards}). The better known approach, sometimes called interpolation (see \cite{bayart2014bohr} or \cite{albuquerque2017holder} for an alternative approach via Hölder's inequality), provides the following estimates (cf. \cite[p.737, Corollary 3.4]{bayart2014bohr}):
\begin{center}
$\begin{array}{|c|c|c|}
\hline
m & \text{best known upper bound for }C_m &  \\ \hline
3 & C_3\leq 1.6818 &  \\ \hline
4 & C_4\leq 1.8877 &  \\ \hline
5 & C_5\leq 2.0586 &  \\ \hline
6 & C_6\leq 2.2064 &  \\ \hline
7 & C_7\leq 2.3376 &  \\ \hline
8 & C_8\leq 2.4562 &  \\ \hline
\end{array}$ \\

\

Table 1
\end{center}

This approach, however, doesn't seem to be optimal, since it appears to lose
information. For a discussion about this, we refer to \cite{pellegrino2017towards}. On the other hand, the lower bounds $2^{1-1/m}$, $m\in\mathbb{N}$, obtained in 2014 (see \cite{diniz2014lower}) were never
improved. In \cite{pellegrino2017towards}, the authors propose a conjecture
that these are in fact the optimal Bohnenblust-Hille constants (see p. 4). A
new approach for the lower bounds, proposed in \cite{cavalcante2017on}, is
capable to provide the truncated Bohnenblust-Hille constants of any order
(although the method is computationally costly). With this new method we are
able to corroborate or refute this conjecture. In addition, it enables us to
elaborate new ones. Through this paper we deal with this approach and its
consequences and issues while we try to implement it computationally.

\section{New approach to the constants}

The new approach introduced in \cite{cavalcante2017on} for the truncated constants consists of finding the finite
set $\text{ext}(B_{\mathcal{L}(^m\mathbb{R}^n)})$ of the extreme points of
the closed unit ball $B_{\mathcal{L}(^m\mathbb{R}^n)}$ in the space ${\mathcal{L}(^m\mathbb{R}^n)}$ of the $m$-linear forms $T:\mathbb{R}^n\times\overset{m}{\cdots} \times\mathbb{R}^n\to\mathbb{R}$. In fact, from (\ref{bhtruncatedinequality}), we have to solve the optimization problem
\begin{equation}  \label{optimization1}
C_m^{(n)}=\sup_{T\in B_{\mathcal{L}(^m\mathbb{R}^n)}}\biggl(\sum_{i_1, ...,
i_m}^n |T(e_{i_1},..., e_{i_m})|^{\frac{2m}{m+1}}\biggr)^{\frac{m+1}{2m}}.
\end{equation}
From the Krein-Milman Theorem, the supremum (\ref{optimization1}) is
attained in one of the extreme points of $B_{\mathcal{L}(^m\mathbb{R}^n)}$,
that is, we reduce the problem of finding the truncated Bohnenblust-Hille
constants to the problem of calculating
\begin{equation}  \label{optimization2}
C_m^{(n)}=\max_{T\in \text{ext}(B_{\mathcal{L}(^m\mathbb{R}^n)})}\biggl(\sum_{i_1, ..., i_m}^n |T(e_{i_1},..., e_{i_m})|^{\frac{2m}{m+1}}\biggr)^{\frac{m+1}{2m}}.
\end{equation}
This way, it is sufficient to know the set $\text{ext}(B_{\mathcal{L}(^m\mathbb{R}^n)})$ and, for this, we can apply the algorithm described in the
next section.

\section{A characterization for the extreme points of $B_{\mathcal{L}(^m\mathbb{R}^n)}$}

From the Theorem 15 of \cite{cavalcante2017on}, we got a full
characterization of the set $\text{ext}(B_{\mathcal{L}(^m\mathbb{R}^n)})$
for all $m,n\in\mathbb{N}$. Now we will summarize the necessary definitions and
state the theorem properly. More details and the proofs can be found within
the reference \cite{cavalcante2017on}.

Given $m,n\in\mathbb{N}$ and $x=(x_1,...,x_m)$, where $x_i\in\mathbb{R}^n$
for $i=1,...,m$, define $I_n\coloneqq \{1,...,n\}$ and, writing the elements
of the cartesian $I_n^m$ in the form $\mathbf{j}=(j_1,...,j_m)$, put
\begin{equation*}
x^{\mathbf{j}}\coloneqq x_1^{(j_1)}\cdots x_m^{(j_m)}
\end{equation*}
and
\begin{equation*}
\omega(x)\coloneqq (x^{\mathbf{j}})_{\mathbf{j}\in I_n^m}\in\mathbb{R}^{n^m},
\end{equation*}
where the $n^m$-tuple $\omega(x)$ is ordered lexicographically. Also, let
\begin{equation}  \label{defVnm}
V_n^m\coloneqq\{\omega(x)\in\mathbb{R}^{n^m}: x\in[\text{ext}(B_{\mathbb{R}^{n}})]^m\}
\end{equation}
and
\begin{equation}  \label{defGnm}
G_n^m\coloneqq \{\text{diag}(v): v\in V_n^m\},
\end{equation}
where $\text{diag}(a_i)_{i\in I_k}$, for $k\in\mathbb{N}$, denotes the $k\times k$ matrix obtained from a list $(a_1,...,a_k)$ putting it in the
diagonal and filling the other entries with 0. Each $m$-linear form $T\in{\mathcal{L}(^m\mathbb{R}^n)}$ can be represented uniquely in the form
\begin{equation}  \label{writeTas0}
T(x)=\sum_{\mathbf{j}}a_{\mathbf{j}}x^{\mathbf{j}},~\text{for all}~x\in[\mathbb{R}^n]^m,
\end{equation}
for some unique real scalars $a_{\mathbf{j}}\in\mathbb{R}$, $\mathbf{j}\in
I_n^m$. Letting $\mathfrak{m}^T\coloneqq (a_{\mathbf{j}})_{\mathbf{j}\in
I_n^m}$, from (\ref{writeTas0}) we have
\begin{equation}  \label{writeTas1}
T(x)=\langle \mathfrak{m}^T, \omega(x) \rangle,~\text{for all}~x\in[\mathbb{R}^n]^m.
\end{equation}
It is easy to see that $\mathfrak{m}^{T\pm L}=\mathfrak{m}^T\pm\mathfrak{m}^L
$ for all $T, L \in{\mathcal{L}(^m\mathbb{R}^n)}$. Writing the formula (\ref{normTformula0}) in the notation (\ref{writeTas1}), we have
\begin{equation}  \label{normTformula1}
\| T\|=\{|\langle \mathfrak{m}^T,v \rangle|\in\mathbb{R}: v\in V_n^m\},~\text{for all}~T\in{\mathcal{L}(^m\mathbb{R}^n)}.
\end{equation}
So, from (\ref{normTformula1}), we have
\begin{equation}  \label{forumaBLmn}
B_{\mathcal{L}(^m\mathbb{R}^n)}=\{T\in{\mathcal{L}(^m\mathbb{R}^n)}:
|\langle \mathfrak{m}^T,v \rangle|\leq1~\text{for all}~v\in V_n^m\}.
\end{equation}

Now, fixing $m,n\in\mathbb{N}$, let us introduce the notations the will lead
us to the characterization of $\text{ext}(B_{\mathcal{L}(^m\mathbb{R}^n)})$.
Denote by $\mathbf{1}$ the vector $(1,\overset{n^m}{...},1)$. We will refer
to it as the \emph{fixed row}. Let $\mathcal{B}$ be the set of all basis of $\mathbb{R}^{n^m}$ contained in $V_n^m$ and containing $\mathbf{1}$, say $\mathcal{=}\{\beta_1,...,\beta_s\}$, for some $s\in\mathbb{N}$. Enumerate
each $\beta_j$, $j\in I_s$, as $\beta_j=\{v_{j,1},...,v_{j,n^m}\}$ so that
the ``position'' of $\mathbf{1}$ is the same for all $\beta_j, j\in I_s$
(say the last one, that is, $v_{j,n^m}=\mathbf{1}$ for all $j\in I_s$).
Define the matrices
\begin{equation*}
H_{\beta_j}\coloneqq \left(
\begin{array}{ccc}
v_{j, 1}^{(1)} & \cdots & v_{j, 1}^{(n^m)} \\
\vdots & \ddots & \vdots \\
v_{j, n^m}^{(1)} & \cdots & v_{j, n^m}^{(n^m)}\end{array}
\right),~\text{for all}~j\in I_s,
\end{equation*}
and the set
\begin{equation*}
\mathcal{A}\coloneqq\{\mathfrak{m}^T\in\mathbb{R}^{n^m}: H_{\beta_j}(\mathfrak{m}^T)^t=e^t, ~\text{for some}~j\in I_s~\text{and}~e\in \text{ext}(B_{\mathbb{R}^{n^m}})\}.
\end{equation*}
The elements of $\mathcal{A}$ are precisely those of the form $\mathfrak{m}^T=(H_{\beta_j}^{-1}e^t)^t$, where $e\in \text{ext}(B_{\mathbb{R}^{n^m}})$
and $j\in I_s$. Finally, define
\begin{equation}  \label{definitionofC}
\mathcal{C}\coloneqq \{\mathfrak{m}^T\cdot g: \mathfrak{m}^T\in\mathcal{A},
T\in B_{\mathcal{L}(^m\mathbb{R}^n)}~\text{and}~g\in G_n^m\}.
\end{equation}
From \cite{cavalcante2017on}, Theorem 15, we conclude that
\begin{equation}  \label{characequality}
\mathcal{C}=\text{ext}(B_{\mathcal{L}(^m\mathbb{R}^n)}).
\end{equation}

\section{The algorithm}

\label{secalgorithm}

As we saw in the last section, we already have an algorithm to calculate the
set $\text{ext}(B_{\mathcal{L}(^m\mathbb{R}^n)})$, which one is composed by
finitely many elementary steps. Before summarize the process, let us reduce
the construction of the set $\mathcal{C}$ defined in (\ref{definitionofC}).

Fix $m,n\in\mathbb{N}$. It is easy to verify that the sets $V_n^m$ defined
in (\ref{defVnm}) and, consequently, $G_n^m$, defined in (\ref{defGnm}), are
both symmetrical. Also, it is easy to see that matrices that are the same,
up to symmetrical rows, leads to the same extreme points, up to symmetry.
Further, since $-\mathbf{1}=(-1,\overset{m}{...},-1)\in G_n^m$, the last
step adds all the missing symmetrical extreme points. So we don't need to
consider the complete set $\mathcal{B}$ of all the bases for $\mathbb{R}^{n^m}$ within $V_n^m$ and containing the fixed row $\mathbf{1}$. Instead,
we remove from $V_n^m$ the symmetrical elements and consider a subset $\overline{V}_n^m$ that has the half of the size of $V_n^m$. Now, instead of $\mathcal{B}$, we take the subset $\mathcal{B}_0$ of the bases for $\mathbb{R}^{n^m}$ within $\overline{V}_n^m$ and containing the fixed row. This
reduction has an enormous impact in the number of matrices considered in the
process. For example, for the planar case (that is, when $n=2$), this
reduction leads to a single base $\beta\in\mathcal{B}_0$, that is, the
process will be made with just one matrix. Other reduction in the process
concerns the set $\text{ext}(B_{\mathbb{R}^{n^m}})$ while calculating $\mathcal{A}$: we only need to proceed with the half of it. The argument is
analogously by symmetry and the fact that $-\mathbf{1}\in G_n^m$, what leads
to a subset $\mathcal{A}_0$ with the half of the size.

Now we are able to describe the process that leads to the set $\mathcal{C}$
defined in (\ref{definitionofC}). Let $m,n\in\mathbb{N}$. The constructive
process consists of performing each of the steps below:
\begin{enumerate}[{Step} 1:]
\item Calculate $I_n^m$, $[\ext(B_{\R^{n}})]^m$ and $\omega(x)$ for all $x\in[\ext(B_{\R^{n}})]^m$ and obtain the set $V_n^m$;
\item Obtain the set $\overline{V}_n^m$ from $V_n^m$ by removing its symmetrical elements;
 \item Consider the set $\mathcal{M}$ of all the invertible matrices with lines in $\overline{V}_n^m$, with the fixed row $\1$ as the last line and removing those that are the same under row switching\footnote{In order to save memory, this step can be made in parts this way: (1) start with the 1-element set $\mathcal{M}_1$ composed by the $1\times n^m$ matrix $\1$; (2) construct the set $\mathcal{M}_2$ of all the $2\times n^m$ matrices by prepending the vectors from $\overline{V}_n^m$ to the matrix from $\mathcal{M}_1$, removing those that are the same by row switching and keeping only those with rank 2; (3) construct the set $\mathcal{M}_3$ of all the $3\times n^m$ matrices by prepending the vectors from $\overline{V}_n^m$ to the matrices from $\mathcal{M}_2$, removing those that are the same by row switching and keeping only those with rank 3; (4) repeat this until you get the set $\mathcal{M}_{n^m}=\mathcal{M}$.};
\item For the half of elements from $e\in\ext(B_{\R^{n^m}})$ (disregarding the symmetric ones) and for each matrix obtained from the last step, calculate $mf^T=e(H_{\beta}^{-1})^t$, what gives us the set $\mathcal{A}_0$;
\item Calculate $\mathcal{C}_0=\mathcal{A}_0\cap B_\Lmn$ using (\ref{forumaBLmn});
\item Finally, obtain $\mathcal{C}$ as the collection of the $n^m$-tuples of the form $\mf^T\cdot g$, with $\mf^T\in\mathcal{C}_0$ and $g\in G_n^m$. From (\ref{characequality}), this set is $\ext(B_\Lmn)$.
\end{enumerate}

\section{Computational implementation}

The process described step-by-step in the last section was implemented by us
in Wolfram Language using the software Mathematica. Two codes were built.
The first one takes as inicial inputs the values of $m$ and $n$ and execute
all the steps listed before in the very same way it was described. We will
refer to this one as the \emph{full program}. The second code, on the other
hand, aims to find extreme points in cases where we can't execute the full
program due to computational limitations. This one, which will be referred
as the \emph{random program}, takes as initial inputs the values of $m$, $n$
and three other ones: (1) the minimum quantity of monomials that the extreme
points to be found should have; (2) the quantity of randomly generated
points within $\text{ext}(B_{\mathbb{R}^{n^m}})$ per test; and (3) the
number of times each matrix randomly generated will be tested. Given these
inputs, the random program execute the following steps:
\begin{enumerate}[{Step} 1:]
\item Execute steps 1 and 2 of the full program;
\item Take a random sample $\{v_1,..., v_{n^m-1}\}$ of $n^m-1$ elements from $\overline{V}_n^m\backslash {\1}$, consider the matrix $H$ whose rows are the vectors from $\{v_1,..., v_{n^m-1}, \1\}$ and calculate its determinant in order to verify if it is invertible;
\item If $H$ is invertible, proceed to the next step, if not, generate another random sample and repeat until the matrix $H$ obtained is invertible;
\item Generate randomly a set of points from $\ext(B_{\R^{n^m}})$ (the quantity decided in the inicial input), calculate $T_{e,H}=(H^{-1}e^t)^t$ with them and collect (using \ref{forumaBLmn}) those that belongs to $B_\Lmn$;
\item If the set collected in the last step is empty, generate new randomly points from $\ext(B_{\R^{n^m}})$ and repeat until some form $T_{e,H}$ within $B_\Lmn$ is found;
\item Within the forms found in the last step, collect those whose quantity of monomials match the minimum required in the inicial input and count this as one test with this matrix\footnote{The emptiness of the set collected in this step is an indicative that the matrix generated is not ``good''. To be sure, we test it more times before generate another one.};
\item Repeat again the steps from 4 forward the number of times decided in the inicial output;
\item If there is some remaining point (what will be extreme points with a interesting quantity of monomials), the algorithm ends, if not, generate another matrix and repeat until something is found.
\end{enumerate}
The \emph{Mathematica notebook} files (.nd) containing these two codes can
be found in the link \href{https://goo.gl/MsdyUV}{https://goo.gl/MsdyUV}. %within the support files attached to this paper.
One observation is that both the notebooks have, at the end, a piece of code aiming to evaluate
the Bohnenblust-Hille constants based in the results obtained with the
programs above. With the full program, the result is the exact truncated
Bohnenblust-Hille constant, while with the random program, the result is a
lower bound. More about this piece of code and the geral output of the
programs will be discussed later in this paper.

\section{About the conjecture}

Fixing $m\in\mathbb{N}$, as we already discussed, the better known lower
bounds for the Bohnen-blust-Hille constants $C_m$ are $2^{1-1/m}$. The
discussion and conjecture that these are the optimal ones first appeared at
\cite{pellegrino2017towards}. There, the authors proved that it is necessary
forms with $4^{m-1}$ monomials, at least, in order to overcome this value.
They key seems to have forms with many monomials and little norm. For
example, for $m=3$, we need forms with more than 16 monomials to overcome
the better known lower bound $2^{1-1/3}$. If this is not the optimal
constant, that is, if $2^{1-1/3}<C_m=\sup \{C_3^{(n)}: n\in\mathbb{N}\}$,
certainly there exists one point in $\text{ext}(B_{\mathcal{L}(^3\mathbb{R}^n)})$, for some $n\in\mathbb{N}$, with more than 16 monomials. However, the
existence of these points remain unknown. On the other hand, if there are no
extreme points in $\text{ext}(B_{\mathcal{L}(^3\mathbb{R}^n)}$, for all $n\in\mathbb{N}$, with more than $16$ monimials, then the conjecture is confirmed
and the optimal Bohnenblust-Hille constant is $2^{1-1/3}$. Obviously this
discussion applies to any value of $m\in\mathbb{N}$ rather than $m=3$.

There are some possibilities with the approach we aim to apply. The first
one and more hopeful is to corroborate the conjecture, what can happens in
many ways. For example: (1) there are no extreme points with more than $4^{m-1}$ monomials; or (2) even if there are some, they don't overcome the
conjectured value. Other possibility with this approach is to find an
extreme point that attain a known upper bound, what is sufficient to find
the optimal constant for the case of $m$.

Speaking about the possibilities with our codes, they allow us to: (1) put
the approach in practice; (2) find extreme points for $B_{\mathcal{L}(^m\mathbb{R}^n)}$ in cases where they were never found; (3) analyse the
behavior of the extreme points in association with the Bohnenblust-Hille
inequality; and (4) investigate the existence of extreme points in $B_{\mathcal{L}(^m\mathbb{R}^n)}$ with more than $4^{m-1}$ monomials.

\section{Investigating the Bohnenblust-Hille constants}

This section aims to describe the last piece of code included in both full
and random programs. Given $m,n\in\mathbb{N}$, the aim of the code is to
calculate the value of the convex function $f:B_{\mathcal{L}(^m\mathbb{R}^n)}\to\mathbb{R}$ given by
\begin{equation}  \label{calculations}
f(T)=\biggl(\sum_{i_1, ..., i_m=1}^n |T(e_{i_1},..., e_{i_m})|^{\frac{2m}{m+1}}\biggr)^{\frac{m+1}{2m}}
\end{equation}
with the extreme points $T$ found within $B_{\mathcal{L}(^m\mathbb{R}^n)}$
and then take the maximum value. With the full program, this process leads
to the truncated constant (\ref{optimization2}), and with the random
program, it leads to a lower bound for it and, consequently, for the optimal
constant $m\in\mathbb{N}$.

As we already mentioned, this last piece of code was included at the end of
both programs. The general outputs are:
\begin{enumerate}[{Out} 1:]
\item The maximum value obtained, say $M\in\R$;
\item Evaluation (true-false) of the the sentence $M\leq2^{1-1/m}$;
\item The quantity of extreme points found;
\item The list of extreme points in matricial form;
\item the list of calculations (\ref{calculations}) made with each extreme points found (exhibited in the respective order in relation to Out 4).
\end{enumerate}

If Out 2 returns False, this means that we refute the conjecture. Searching
in Out 5 the user can find the extreme point that leads the maximum result,
pick its position and look at the extreme point in Out 4 at the same
position. This is the interesting form found though the program. We can also
easily count the number of monomials in each extreme point found and see if
something interesting appear.

\section{Some results and new conjectures}

Since, for each $m\in\mathbb{N}$, we are looking for extreme points with
more than $4^{m-1}$ monomials and in ${\mathcal{L}(^m\mathbb{R}^n)}$ the
forms have $n^m$ at maximum, it is interesting to fix $m$ and take $n$ such
that $n^m>4^{m-1}$, that is, $n>\sqrt[m]{4^{m-1}}$. Since we aim to
investigate the Bohnenblust-Hille constants and since for $m=2$ we recover
the 4/3-Littlewood inequality, for which we already know that $C_m=\sqrt{2}$, it is interesting for us to run the case $m\geq3$. So the smallest
interesting cases are $n\geq3$.

With the computational resources available for us, we couldn't run the full
program for $m\geq3$ and $n\geq 3$ (although for smaller values of $m$ or $n$
it was possible). The issues are discussed in the last section of this
paper. We could, however, adjust the random program to fit our available
resources and find some interesting results.

\begin{figure}[!htb]
	%\vspace*{0,2cm}
    \centering
    \caption{Graph of $n\mapsto f(T_n)$.}
    \includegraphics[width=13cm]{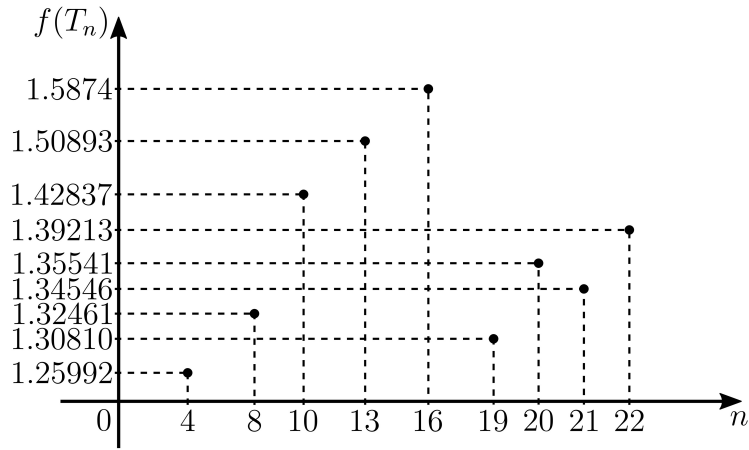}
    \label{figmonomios}
    \medskip
    \small
\end{figure}

Configuring the random program to test each matrix 50 times and generate
10000 points from $\text{ext}(B_{\mathbb{R}^{3^3}})$ per evaluation, we
could find some interesting points. While searching for extreme with at
least 16 monomials, the program find a lot of them easily. For example, the
form
\begin{equation*}
T_{16}=\frac{1}{4}(-1,-1,0,-1,-1,0,0,0,0,-1,0,-1,0,1,1,-1,1,0,0,-1,1,1,0,-1,-1,1,0)
\end{equation*}
is one of the extreme points found. All of them attained the conjectured ($2^{1-1/3}\approx1.5874$) value when we calculated (\ref{calculations}).
Although doubting the existence of extreme points in $B_{\mathcal{L}(^3\mathbb{R}^3)}$ with more than 16 monomials, we ran the code searching
extremes with at least 17 monomials. After about 6 hours, for our surprise,
there were found two points. The first, with 19 monomials, is
\begin{equation*}
T_{19}=-\frac{1}{8}(1, 2, -1, 0, 0, 0, -1, 2, -1, 1, 1, 4, 1, 0, -1, 0, 1,
-3, 0, 1, 1, -1, 0, 1, 1, 1, 0).
\end{equation*}
The second, with 20 monomials, is
\begin{equation*}
T_{20}=-\frac{1}{8}(1, 0, -1, 0, 0, 0, 1, 0,- 1, 1, -1, -2, 1, 2, 1, 0, 1,
1, 2, 1, -1, 3, -2, 3, -1, -1, 0).
\end{equation*}
Another interesting fact about this discovery, rather than the existence,
was the values of (\ref{calculations}) for these two monomials. We got
\begin{equation*}
f(T_{19})\approx 1.30810~~\text{e}~~f(T_{20})\approx 1.35541.
\end{equation*}
That is, $f(T_{19})<f(T_{20})<f(T_{16})$. The ``good'' matrix that generated
these points is \scriptsize
\begin{equation*}
\left(\begin{array}{rrrrrrrrrrrrrrrrrrrrrrrrrrr}
\!\! 1 \!\! & \!\! 1 \!\! & \!\! 1 \!\! & \!\! 1 \!\! & \!\! 1 \!\! & \!\! 1
\!\! & \!\! 1 \!\! & \!\! 1 \!\! & \!\! 1 \!\! & \!\! 1 \!\! & \!\! 1 \!\! &
\!\! 1 \!\! & \!\! 1 \!\! & \!\! 1 \!\! & \!\! 1 \!\! & \!\! 1 \!\! & \!\! 1
\!\! & \!\! 1 \!\! & \!\! 1 \!\! & \!\! 1 \!\! & \!\! 1 \!\! & \!\! 1 \!\! &
\!\! 1 \!\! & \!\! 1 \!\! & \!\! 1 \!\! & \!\! 1 \!\! & \!\! 1 \\
\!\! 1 \!\! & \!\! -1 \!\! & \!\! -1 \!\! & \!\! -1 \!\! & \!\! 1 \!\! &
\!\! 1 \!\! & \!\! -1 \!\! & \!\! 1 \!\! & \!\! 1 \!\! & \!\! 1 \!\! & \!\!
-1 \!\! & \!\! -1 \!\! & \!\! -1 \!\! & \!\! 1 \!\! & \!\! 1 \!\! & \!\! -1
\!\! & \!\! 1 \!\! & \!\! 1 \!\! & \!\! 1 \!\! & \!\! -1 \!\! & \!\! -1 \!\!
& \!\! -1 \!\! & \!\! 1 \!\! & \!\! 1 \!\! & \!\! -1 \!\! & \!\! 1 \!\! &
\!\! 1 \\
\!\! 1 \!\! & \!\! 1 \!\! & \!\! -1 \!\! & \!\! -1 \!\! & \!\! -1 \!\! &
\!\! 1 \!\! & \!\! -1 \!\! & \!\! -1 \!\! & \!\! 1 \!\! & \!\! -1 \!\! &
\!\! -1 \!\! & \!\! 1 \!\! & \!\! 1 \!\! & \!\! 1 \!\! & \!\! -1 \!\! & \!\!
1 \!\! & \!\! 1 \!\! & \!\! -1 \!\! & \!\! 1 \!\! & \!\! 1 \!\! & \!\! -1
\!\! & \!\! -1 \!\! & \!\! -1 \!\! & \!\! 1 \!\! & \!\! -1 \!\! & \!\! -1
\!\! & \!\! 1 \\
\!\! 1 \!\! & \!\! -1 \!\! & \!\! 1 \!\! & \!\! -1 \!\! & \!\! 1 \!\! & \!\!
-1 \!\! & \!\! 1 \!\! & \!\! -1 \!\! & \!\! 1 \!\! & \!\! -1 \!\! & \!\! 1
\!\! & \!\! -1 \!\! & \!\! 1 \!\! & \!\! -1 \!\! & \!\! 1 \!\! & \!\! -1 \!\!
& \!\! 1 \!\! & \!\! -1 \!\! & \!\! -1 \!\! & \!\! 1 \!\! & \!\! -1 \!\! &
\!\! 1 \!\! & \!\! -1 \!\! & \!\! 1 \!\! & \!\! -1 \!\! & \!\! 1 \!\! & \!\!
-1 \\
\!\! 1 \!\! & \!\! 1 \!\! & \!\! 1 \!\! & \!\! 1 \!\! & \!\! 1 \!\! & \!\! 1
\!\! & \!\! -1 \!\! & \!\! -1 \!\! & \!\! -1 \!\! & \!\! -1 \!\! & \!\! -1
\!\! & \!\! -1 \!\! & \!\! -1 \!\! & \!\! -1 \!\! & \!\! -1 \!\! & \!\! 1
\!\! & \!\! 1 \!\! & \!\! 1 \!\! & \!\! 1 \!\! & \!\! 1 \!\! & \!\! 1 \!\! &
\!\! 1 \!\! & \!\! 1 \!\! & \!\! 1 \!\! & \!\! -1 \!\! & \!\! -1 \!\! & \!\!
-1 \\
\!\! 1 \!\! & \!\! -1 \!\! & \!\! 1 \!\! & \!\! -1 \!\! & \!\! 1 \!\! & \!\!
-1 \!\! & \!\! -1 \!\! & \!\! 1 \!\! & \!\! -1 \!\! & \!\! -1 \!\! & \!\! 1
\!\! & \!\! -1 \!\! & \!\! 1 \!\! & \!\! -1 \!\! & \!\! 1 \!\! & \!\! 1 \!\!
& \!\! -1 \!\! & \!\! 1 \!\! & \!\! -1 \!\! & \!\! 1 \!\! & \!\! -1 \!\! &
\!\! 1 \!\! & \!\! -1 \!\! & \!\! 1 \!\! & \!\! 1 \!\! & \!\! -1 \!\! & \!\!
1 \\
\!\! 1 \!\! & \!\! -1 \!\! & \!\! 1 \!\! & \!\! -1 \!\! & \!\! 1 \!\! & \!\!
-1 \!\! & \!\! 1 \!\! & \!\! -1 \!\! & \!\! 1 \!\! & \!\! -1 \!\! & \!\! 1
\!\! & \!\! -1 \!\! & \!\! 1 \!\! & \!\! -1 \!\! & \!\! 1 \!\! & \!\! -1 \!\!
& \!\! 1 \!\! & \!\! -1 \!\! & \!\! 1 \!\! & \!\! -1 \!\! & \!\! 1 \!\! &
\!\! -1 \!\! & \!\! 1 \!\! & \!\! -1 \!\! & \!\! 1 \!\! & \!\! -1 \!\! &
\!\! 1 \\
\!\! 1 \!\! & \!\! -1 \!\! & \!\! 1 \!\! & \!\! -1 \!\! & \!\! 1 \!\! & \!\!
-1 \!\! & \!\! -1 \!\! & \!\! 1 \!\! & \!\! -1 \!\! & \!\! -1 \!\! & \!\! 1
\!\! & \!\! -1 \!\! & \!\! 1 \!\! & \!\! -1 \!\! & \!\! 1 \!\! & \!\! 1 \!\!
& \!\! -1 \!\! & \!\! 1 \!\! & \!\! 1 \!\! & \!\! -1 \!\! & \!\! 1 \!\! &
\!\! -1 \!\! & \!\! 1 \!\! & \!\! -1 \!\! & \!\! -1 \!\! & \!\! 1 \!\! &
\!\! -1 \\
\!\! 1 \!\! & \!\! 1 \!\! & \!\! -1 \!\! & \!\! -1 \!\! & \!\! -1 \!\! &
\!\! 1 \!\! & \!\! -1 \!\! & \!\! -1 \!\! & \!\! 1 \!\! & \!\! -1 \!\! &
\!\! -1 \!\! & \!\! 1 \!\! & \!\! 1 \!\! & \!\! 1 \!\! & \!\! -1 \!\! & \!\!
1 \!\! & \!\! 1 \!\! & \!\! -1 \!\! & \!\! -1 \!\! & \!\! -1 \!\! & \!\! 1
\!\! & \!\! 1 \!\! & \!\! 1 \!\! & \!\! -1 \!\! & \!\! 1 \!\! & \!\! 1 \!\!
& \!\! -1 \\
\!\! 1 \!\! & \!\! -1 \!\! & \!\! 1 \!\! & \!\! 1 \!\! & \!\! -1 \!\! & \!\!
1 \!\! & \!\! -1 \!\! & \!\! 1 \!\! & \!\! -1 \!\! & \!\! -1 \!\! & \!\! 1
\!\! & \!\! -1 \!\! & \!\! -1 \!\! & \!\! 1 \!\! & \!\! -1 \!\! & \!\! 1 \!\!
& \!\! -1 \!\! & \!\! 1 \!\! & \!\! -1 \!\! & \!\! 1 \!\! & \!\! -1 \!\! &
\!\! -1 \!\! & \!\! 1 \!\! & \!\! -1 \!\! & \!\! 1 \!\! & \!\! -1 \!\! &
\!\! 1 \\
\!\! 1 \!\! & \!\! 1 \!\! & \!\! 1 \!\! & \!\! -1 \!\! & \!\! -1 \!\! & \!\!
-1 \!\! & \!\! -1 \!\! & \!\! -1 \!\! & \!\! -1 \!\! & \!\! 1 \!\! & \!\! 1
\!\! & \!\! 1 \!\! & \!\! -1 \!\! & \!\! -1 \!\! & \!\! -1 \!\! & \!\! -1
\!\! & \!\! -1 \!\! & \!\! -1 \!\! & \!\! 1 \!\! & \!\! 1 \!\! & \!\! 1 \!\!
& \!\! -1 \!\! & \!\! -1 \!\! & \!\! -1 \!\! & \!\! -1 \!\! & \!\! -1 \!\! &
\!\! -1 \\
\!\! 1 \!\! & \!\! -1 \!\! & \!\! -1 \!\! & \!\! -1 \!\! & \!\! 1 \!\! &
\!\! 1 \!\! & \!\! -1 \!\! & \!\! 1 \!\! & \!\! 1 \!\! & \!\! 1 \!\! & \!\!
-1 \!\! & \!\! -1 \!\! & \!\! -1 \!\! & \!\! 1 \!\! & \!\! 1 \!\! & \!\! -1
\!\! & \!\! 1 \!\! & \!\! 1 \!\! & \!\! -1 \!\! & \!\! 1 \!\! & \!\! 1 \!\!
& \!\! 1 \!\! & \!\! -1 \!\! & \!\! -1 \!\! & \!\! 1 \!\! & \!\! -1 \!\! &
\!\! -1 \\
\!\! 1 \!\! & \!\! 1 \!\! & \!\! 1 \!\! & \!\! 1 \!\! & \!\! 1 \!\! & \!\! 1
\!\! & \!\! -1 \!\! & \!\! -1 \!\! & \!\! -1 \!\! & \!\! -1 \!\! & \!\! -1
\!\! & \!\! -1 \!\! & \!\! -1 \!\! & \!\! -1 \!\! & \!\! -1 \!\! & \!\! 1
\!\! & \!\! 1 \!\! & \!\! 1 \!\! & \!\! -1 \!\! & \!\! -1 \!\! & \!\! -1 \!\!
& \!\! -1 \!\! & \!\! -1 \!\! & \!\! -1 \!\! & \!\! 1 \!\! & \!\! 1 \!\! &
\!\! 1 \\
\!\! 1 \!\! & \!\! 1 \!\! & \!\! -1 \!\! & \!\! 1 \!\! & \!\! 1 \!\! & \!\!
-1 \!\! & \!\! -1 \!\! & \!\! -1 \!\! & \!\! 1 \!\! & \!\! 1 \!\! & \!\! 1
\!\! & \!\! -1 \!\! & \!\! 1 \!\! & \!\! 1 \!\! & \!\! -1 \!\! & \!\! -1 \!\!
& \!\! -1 \!\! & \!\! 1 \!\! & \!\! 1 \!\! & \!\! 1 \!\! & \!\! -1 \!\! &
\!\! 1 \!\! & \!\! 1 \!\! & \!\! -1 \!\! & \!\! -1 \!\! & \!\! -1 \!\! &
\!\! 1 \\
\!\! 1 \!\! & \!\! -1 \!\! & \!\! -1 \!\! & \!\! 1 \!\! & \!\! -1 \!\! &
\!\! -1 \!\! & \!\! 1 \!\! & \!\! -1 \!\! & \!\! -1 \!\! & \!\! -1 \!\! &
\!\! 1 \!\! & \!\! 1 \!\! & \!\! -1 \!\! & \!\! 1 \!\! & \!\! 1 \!\! & \!\!
-1 \!\! & \!\! 1 \!\! & \!\! 1 \!\! & \!\! 1 \!\! & \!\! -1 \!\! & \!\! -1
\!\! & \!\! 1 \!\! & \!\! -1 \!\! & \!\! -1 \!\! & \!\! 1 \!\! & \!\! -1 \!\!
& \!\! -1 \\
\!\! 1 \!\! & \!\! -1 \!\! & \!\! 1 \!\! & \!\! -1 \!\! & \!\! 1 \!\! & \!\!
-1 \!\! & \!\! 1 \!\! & \!\! -1 \!\! & \!\! 1 \!\! & \!\! 1 \!\! & \!\! -1
\!\! & \!\! 1 \!\! & \!\! -1 \!\! & \!\! 1 \!\! & \!\! -1 \!\! & \!\! 1 \!\!
& \!\! -1 \!\! & \!\! 1 \!\! & \!\! 1 \!\! & \!\! -1 \!\! & \!\! 1 \!\! &
\!\! -1 \!\! & \!\! 1 \!\! & \!\! -1 \!\! & \!\! 1 \!\! & \!\! -1 \!\! &
\!\! 1 \\
\!\! 1 \!\! & \!\! 1 \!\! & \!\! -1 \!\! & \!\! 1 \!\! & \!\! 1 \!\! & \!\!
-1 \!\! & \!\! 1 \!\! & \!\! 1 \!\! & \!\! -1 \!\! & \!\! 1 \!\! & \!\! 1
\!\! & \!\! -1 \!\! & \!\! 1 \!\! & \!\! 1 \!\! & \!\! -1 \!\! & \!\! 1 \!\!
& \!\! 1 \!\! & \!\! -1 \!\! & \!\! -1 \!\! & \!\! -1 \!\! & \!\! 1 \!\! &
\!\! -1 \!\! & \!\! -1 \!\! & \!\! 1 \!\! & \!\! -1 \!\! & \!\! -1 \!\! &
\!\! 1 \\
\!\! 1 \!\! & \!\! 1 \!\! & \!\! -1 \!\! & \!\! 1 \!\! & \!\! 1 \!\! & \!\!
-1 \!\! & \!\! 1 \!\! & \!\! 1 \!\! & \!\! -1 \!\! & \!\! 1 \!\! & \!\! 1
\!\! & \!\! -1 \!\! & \!\! 1 \!\! & \!\! 1 \!\! & \!\! -1 \!\! & \!\! 1 \!\!
& \!\! 1 \!\! & \!\! -1 \!\! & \!\! 1 \!\! & \!\! 1 \!\! & \!\! -1 \!\! &
\!\! 1 \!\! & \!\! 1 \!\! & \!\! -1 \!\! & \!\! 1 \!\! & \!\! 1 \!\! & \!\!
-1 \\
\!\! 1 \!\! & \!\! -1 \!\! & \!\! 1 \!\! & \!\! 1 \!\! & \!\! -1 \!\! & \!\!
1 \!\! & \!\! -1 \!\! & \!\! 1 \!\! & \!\! -1 \!\! & \!\! -1 \!\! & \!\! 1
\!\! & \!\! -1 \!\! & \!\! -1 \!\! & \!\! 1 \!\! & \!\! -1 \!\! & \!\! 1 \!\!
& \!\! -1 \!\! & \!\! 1 \!\! & \!\! 1 \!\! & \!\! -1 \!\! & \!\! 1 \!\! &
\!\! 1 \!\! & \!\! -1 \!\! & \!\! 1 \!\! & \!\! -1 \!\! & \!\! 1 \!\! & \!\!
-1 \\
\!\! 1 \!\! & \!\! -1 \!\! & \!\! -1 \!\! & \!\! 1 \!\! & \!\! -1 \!\! &
\!\! -1 \!\! & \!\! 1 \!\! & \!\! -1 \!\! & \!\! -1 \!\! & \!\! 1 \!\! &
\!\! -1 \!\! & \!\! -1 \!\! & \!\! 1 \!\! & \!\! -1 \!\! & \!\! -1 \!\! &
\!\! 1 \!\! & \!\! -1 \!\! & \!\! -1 \!\! & \!\! 1 \!\! & \!\! -1 \!\! &
\!\! -1 \!\! & \!\! 1 \!\! & \!\! -1 \!\! & \!\! -1 \!\! & \!\! 1 \!\! &
\!\! -1 \!\! & \!\! -1 \\
\!\! 1 \!\! & \!\! -1 \!\! & \!\! 1 \!\! & \!\! 1 \!\! & \!\! -1 \!\! & \!\!
1 \!\! & \!\! -1 \!\! & \!\! 1 \!\! & \!\! -1 \!\! & \!\! 1 \!\! & \!\! -1
\!\! & \!\! 1 \!\! & \!\! 1 \!\! & \!\! -1 \!\! & \!\! 1 \!\! & \!\! -1 \!\!
& \!\! 1 \!\! & \!\! -1 \!\! & \!\! -1 \!\! & \!\! 1 \!\! & \!\! -1 \!\! &
\!\! -1 \!\! & \!\! 1 \!\! & \!\! -1 \!\! & \!\! 1 \!\! & \!\! -1 \!\! &
\!\! 1 \\
\!\! 1 \!\! & \!\! 1 \!\! & \!\! -1 \!\! & \!\! -1 \!\! & \!\! -1 \!\! &
\!\! 1 \!\! & \!\! 1 \!\! & \!\! 1 \!\! & \!\! -1 \!\! & \!\! -1 \!\! & \!\!
-1 \!\! & \!\! 1 \!\! & \!\! 1 \!\! & \!\! 1 \!\! & \!\! -1 \!\! & \!\! -1
\!\! & \!\! -1 \!\! & \!\! 1 \!\! & \!\! 1 \!\! & \!\! 1 \!\! & \!\! -1 \!\!
& \!\! -1 \!\! & \!\! -1 \!\! & \!\! 1 \!\! & \!\! 1 \!\! & \!\! 1 \!\! &
\!\! -1 \\
\!\! 1 \!\! & \!\! -1 \!\! & \!\! -1 \!\! & \!\! -1 \!\! & \!\! 1 \!\! &
\!\! 1 \!\! & \!\! 1 \!\! & \!\! -1 \!\! & \!\! -1 \!\! & \!\! 1 \!\! & \!\!
-1 \!\! & \!\! -1 \!\! & \!\! -1 \!\! & \!\! 1 \!\! & \!\! 1 \!\! & \!\! 1
\!\! & \!\! -1 \!\! & \!\! -1 \!\! & \!\! -1 \!\! & \!\! 1 \!\! & \!\! 1 \!\!
& \!\! 1 \!\! & \!\! -1 \!\! & \!\! -1 \!\! & \!\! -1 \!\! & \!\! 1 \!\! &
\!\! 1 \\
\!\! 1 \!\! & \!\! 1 \!\! & \!\! -1 \!\! & \!\! 1 \!\! & \!\! 1 \!\! & \!\!
-1 \!\! & \!\! -1 \!\! & \!\! -1 \!\! & \!\! 1 \!\! & \!\! -1 \!\! & \!\! -1
\!\! & \!\! 1 \!\! & \!\! -1 \!\! & \!\! -1 \!\! & \!\! 1 \!\! & \!\! 1 \!\!
& \!\! 1 \!\! & \!\! -1 \!\! & \!\! -1 \!\! & \!\! -1 \!\! & \!\! 1 \!\! &
\!\! -1 \!\! & \!\! -1 \!\! & \!\! 1 \!\! & \!\! 1 \!\! & \!\! 1 \!\! & \!\!
-1 \\
\!\! 1 \!\! & \!\! 1 \!\! & \!\! 1 \!\! & \!\! 1 \!\! & \!\! 1 \!\! & \!\! 1
\!\! & \!\! -1 \!\! & \!\! -1 \!\! & \!\! -1 \!\! & \!\! 1 \!\! & \!\! 1 \!\!
& \!\! 1 \!\! & \!\! 1 \!\! & \!\! 1 \!\! & \!\! 1 \!\! & \!\! -1 \!\! &
\!\! -1 \!\! & \!\! -1 \!\! & \!\! -1 \!\! & \!\! -1 \!\! & \!\! -1 \!\! &
\!\! -1 \!\! & \!\! -1 \!\! & \!\! -1 \!\! & \!\! 1 \!\! & \!\! 1 \!\! &
\!\! 1 \\
\!\! 1 \!\! & \!\! -1 \!\! & \!\! -1 \!\! & \!\! -1 \!\! & \!\! 1 \!\! &
\!\! 1 \!\! & \!\! -1 \!\! & \!\! 1 \!\! & \!\! 1 \!\! & \!\! -1 \!\! & \!\!
1 \!\! & \!\! 1 \!\! & \!\! 1 \!\! & \!\! -1 \!\! & \!\! -1 \!\! & \!\! 1
\!\! & \!\! -1 \!\! & \!\! -1 \!\! & \!\! 1 \!\! & \!\! -1 \!\! & \!\! -1
\!\! & \!\! -1 \!\! & \!\! 1 \!\! & \!\! 1 \!\! & \!\! -1 \!\! & \!\! 1 \!\!
& \!\! 1 \\
\!\! 1 \!\! & \!\! 1 \!\! & \!\! 1 \!\! & \!\! 1 \!\! & \!\! 1 \!\! & \!\! 1
\!\! & \!\! 1 \!\! & \!\! 1 \!\! & \!\! 1 \!\! & \!\! -1 \!\! & \!\! -1 \!\!
& \!\! -1 \!\! & \!\! -1 \!\! & \!\! -1 \!\! & \!\! -1 \!\! & \!\! -1 \!\! &
\!\! -1 \!\! & \!\! -1 \!\! & \!\! -1 \!\! & \!\! -1 \!\! & \!\! -1 \!\! &
\!\! -1 \!\! & \!\! -1 \!\! & \!\! -1 \!\! & \!\! -1 \!\! & \!\! -1 \!\! &
\!\! -1 \\
&  &  &  &  &  &  &  &  &  &  &  &  &  &  &  &  &  &  &  &  &  &  &  &  &  &
\end{array}\right).
\end{equation*}
\normalsize
Believing that this matrix would lead to more interesting
results, we ran a modified version of the full program to calculate all the
extreme points that this matrix generate. There were found 55720 nontrivial
(that is, different than the canonical vectors) extreme points: 2376 with 4
monomials, 3456 with 8 monomials, 18144 with 10 monomials, 20736 with 13
monomials, 6912 with 16 monomials, 768 with 19 monomials, 1536 with 20
monomials, 256 with 21 monomials and 1536 with 22 monomials. Some example
are:
\begin{gather*}
T_{4}=\frac{1}{2}(1,1,0,0,0,0,0,0,0,0,0,0,0,0,0,-1,1,0,0,0,0,0,0,0,0,0,0), \\
T_{8}=\frac{1}{4}(3,0,1,-1,0,1,0,0,0,0,0,0,0,0,0,0,0,0,1,0,-1,1,0,-1,0,0,0)
\\
T_{10}=\frac{1}{4}(-2,-2,0,0,0,0,0,0,0,0,0,0,1,-1,0,-1,1,0,0,0,0,-1,1,0,-1,1,0) \\
T_{13}=\frac{1}{4}(2,1,1,1,0,-1,-1,1,0,0,0,0,0,0,0,0,0,0,0,1,-1,-1,0,1,-1,1,0) \\
T_{16}=\frac{1}{4}(-1,-1,0,0,-1,1,1,0,1,1,1,0,-1,0,1,0,1,1,0,0,0,-1,1,0,1,-1,0) \\
T_{21}=\frac{1}{8}(1,0,1,0,2,2,1,-2,-1,1,1,0,-1,0,-1,0,1,-1,0,-1,1,1,-2,-1,1,1,-4)~~\text{e}
\end{gather*}
\begin{equation*}
T_{22}=\frac{1}{8}(-3,0,-1,0,2,-2,1,2,1,1,1,0,1,0,-1,0,1,1,2,-1,1,-1,2,-1,-1,1,2).
\end{equation*}
With them, we obtained the following results:
\begin{gather*}
f(T_{4})\approx1.25992,~~f(T_{8})\approx1.32461,~~f(T_{10})\approx1.42837,~~f(T_{13})\approx1.50893, \\
f(T_{16})\approx1.5874,~~f(T_{21})\approx1.34546~~\text{e}~~f(T_{22})\approx1.39213.
\end{gather*}
Our discoveries are resumed in the Figure \ref{figmonomios}.

As we saw above, we couldn't refute the conjecture arisen in
\cite{pellegrino2017towards}, but corroborate it. Moreover, some interesting
facts that we found can suggest new conjectures. As we can see in the Figure \ref{figmonomios}, the extreme points with 16 monomials are ``special'' in
the sense that the value of (\ref{calculations}) were bigger with them than
with any other extreme points, even those with more monomials. So there are
in fact extreme points that have more than 16 monomials for the case $m=n=3$. An interesting thing that we can notice, however, is that the extreme that
led to the best value of (\ref{calculations}) have the same coefficient in
every monomial. They are, in essence, the same kind of form found in \cite{diniz2014lower}, that is, with coefficient $1$ or $-1$ in each parcel. This
way, fixing $m,n\in\mathbb{N}$, two ideas come in mind: (1) maybe the
extreme points of $B_{\mathcal{L}(^m\mathbb{R}^n)}$ with more than $4^{m-1}$
monomials don't have the same coefficient for each of its monomials; and (2)
maybe forms whose monomials don't have the same coefficient can't attain the
best known lower bound $2^{1-1/m}$. These two open problems could be
conceived as conjectures. If they were proved, then we would conclude that
the conjecture arisen in \cite{pellegrino2017towards} is true, that is, $C_m=2^{1-1/m}$.

\section{\protect\normalsize The issues and future perspectives}

In this last section we discuss some issues with the actual
code and the reason why we couldn't run the programs for larger values.

There are two big sets that are dealt with in the process
described in the Section \ref{secalgorithm}. The first one is $\text{ext}(B_{\mathbb{R}^{n^m}})$, $m,n\in\mathbb{N}$. For large values of $n$ and $m$,
this set is too big for any computer to hold. However, as we saw in (\ref{extBRnm}) and using the function \verb|IntegerDigits| from the Wolfram
Language (see the code attached), we have a simple solution to calculate
this set one-element-by-one without holding it entirely at once. This way,
the memory used to store $\text{ext}(B_{\mathbb{R}^{n^m}})$ is not a
problem.

On the other hand, there is a second and problematic set that
we have to deal with in the process, which is $\mathcal{B}$, that is, the
collection of all bases within $\overline{V}_n^m$ containing the fixed line $\mathbf{1}$. The way we constructed it in the Section (\ref{secalgorithm})
is far from be the best. This process cannot be parallelized and the memory
consumption is very high. So we need to find a way to avoid this problem,
maybe calculating the elements of the set $\mathcal{B}$ one-by-one. For
this, we need to better understand the set $V_n^m$ and see what we can do to
generate the matrices. Solving this, the only limitation will be the
computing power of the machine, so we can take full advance of a cluster
while running the full program. In such manner it will be passible to
investigate further the behavior of the truncated constants, while we track
the knowledge of the best ones.

%\section*{Acknowledgments}

%\appendix

%\section{Source codes}

%The elaborated codes previously described (``Full program'' and ``Random program'') are available (in \emph{Mathematica notebook} files ``.nd'') at \href{https://goo.gl/MsdyUV}{https://goo.gl/MsdyUV}.

\bigskip

\bibliographystyle{plain}
%\bibliography{refbase}

\

\end{document}